%% file: counting_abelian_squares.tex
\documentclass[english,11pt,prl]{revtex4-2}
\usepackage[T1]{fontenc}
\usepackage[latin9]{inputenc}
\usepackage[letterpaper]{geometry}
\usepackage{amsmath}
\usepackage{graphicx}

\makeatletter

\providecommand{\tabularnewline}{\\}

\makeatother

\usepackage{babel}
\begin{document}
\input{macros.tex}

\title{Counting Abelian Squares More Efficiently}
\author{Ryan Bennink}
\address{Quantum Computational Science Group, Oak Ridge National Laboratory}
\begin{abstract}
I present a recursive formula for calculating the number of abelian
squares of length $n+n$ over an alphabet of size $d$. The presented
formula is similar to a previously known formula but has substantially
lower complexity when $d\gg n$.
\end{abstract}
\maketitle

\section{Introduction}

An abelian square is a word whose first half is an anagram of its
second half, for example $\mathtt{intestines}=\mathtt{intes}\cdot\mathtt{tines}$
or $\mathtt{bonbon}=\mathtt{bon}\cdot\mathtt{bon}$. Abelian squares
are fundamentally interesting combinatoric structures \cite{Erdos1957,Keranen1992,Iliopoulos1997,Carpi1998,Cassaigne2011,Huova2012,Crochemore2013}
that can arise in a variety of contexts in applied mathematics. The
work reported here was motivated by a problem in quantum computing.
As discussed in a related manuscript \cite{Bennink2022b}, the expressiveness
of a certain class of variational quantum circuits can be related
to the number of abelian squares over a certain alphabet. However,
due to fact that the alphabet in this case is exponentially large,
existing formulas for evaluating the number of abelian squares were
found to be impractical.

In this letter I present a recursive formula for calculating the number
of abelian squares of length $n+n$ over an alphabet of size $d$
that is efficient even when $d$ is very large. I first review the
problem of counting abelian squares and an existing recursive formula
\cite{Richmond2009} whose cost is $O(n^{2}d)$. Then I derive a new
recursive formula whose cost is only $O(n^{2}\min(n,d))$, a substantial
improvement when $d\gg n$. I furthermore give a constructive interpretation
of the formula.

\section{Background}

Let $f_{d}(n)$ denote the number of abelian squares of length $n+n$
over an alphabet of $d$ symbols. Trivially, $f_{1}(n)=1$ for all
$n$ and $f_{d}(0)=1$ for all $d$. It is also not difficult to see
that $f_{d}(1)=d$. To determine $f_{d}(n)$ for arbitrary $d$ and
$n$, we define the \emph{signature} (sometimes called the Parikh
vector) of a word $w\in\{a_{1},\ldots,a_{d}\}^{*}$ as $(m_{1},\ldots,m_{d})$
where $m_{i}$ is the number of times the symbol $a_{i}$ appears
in $w$. Note that two words are anagrams if and only if they have
the same signature. Thus the number of abelian squares is the number
of pairs $(x,y)$ such that $x$ and $y$ have the same signature.
The number of words with a particular signature $(m_{1},\ldots,m_{d})$
is given by the multinomial coefficient
\begin{align}
\binom{m_{1}+\cdots+m_{d}}{m_{1},\ldots,m_{d}} & =\frac{(m_{1}+\cdots+m_{d})!}{m_{1}!\cdots m_{d}!}.
\end{align}
The number of ways to choose a pair of words, each with signature
$(m_{1},\ldots,m_{d})$, is just the square of this quantity. Therefore
the number of abelian squares of length $n+n$ is
\begin{align}
f_{d}(n) & =\sum_{0\le m_{1}+\cdots+m_{d}\le n}\binom{n}{m_{1},\ldots,m_{d}}^{2}.\label{eq: basic formula for AS count}
\end{align}
The first few values of $f_{d}(n)$ are shown in Table \ref{tab: first few values}.

\begin{table}
\begin{centering}
\begin{tabular}{|c||c|c|c|c|c|c|c|c|}
\hline 
$d$\textbackslash$n$ & ~0~ & ~1~ & 2 & 3 & 4 & 5 & 6 & 7\tabularnewline
\hline 
\hline 
1 & 1 & 1 & 1 & 1 & 1 & 1 & 1 & 1\tabularnewline
\hline 
2 & 1 & 2 & 6 & 20 & 70 & 252 & 924 & 3432\tabularnewline
\hline 
3 & 1 & 3 & 15 & 93 & 639 & 4653 & 35169 & 272835\tabularnewline
\hline 
4 & 1 & 4 & 28 & 256 & 2716 & 31504 & 387136 & 4951552\tabularnewline
\hline 
5 & 1 & 5 & 45 & 545 & 7885 & 127905 & 2241225 & 41467725\tabularnewline
\hline 
6 & 1 & 6 & 66 & 996 & 18306 & 384156 & 8848236 & 218040696\tabularnewline
\hline 
\end{tabular}
\par\end{centering}
\caption{\label{tab: first few values}Number of abelian squares of length
$n+n$ over an alphabet of size $d$ \cite{Richmond2009}.}
\end{table}

Eq.~(\ref{eq: basic formula for AS count}) is not easy to evaluate
when $n$ is large, as the number of signatures grows combinatorially
in $d$ and $n$. Richmond and Shallit \cite{Richmond2009} derived
a recursive formula using a simple constructive argument: To create
size $(n,n)$ abelian word pair $(x,y)$ over alphabet $\{a_{1},\ldots,a_{d}\}$,
one can first choose the number $i\in\{0,\ldots,n\}$ of occurrences
of $a_{d}$ in each word. There are $\tbinom{n}{i}$ ways to distribute
these occurrences in each word. Then there are $f_{d-1}(n-i)$ ways
to create an abelian pair over $\{a_{1},\ldots,a_{d-1}\}$ for the
remaining $n-i$ symbols in each word. Setting $k=n-i$ and summing
over the choice of $k$ yields
\begin{align}
f_{d}(n) & =\sum_{k=0}^{n}\binom{n}{k}^{2}f_{d-1}(k).\label{eq: Richmond}
\end{align}
Using this formula, $f_{d}(n)$ can be obtained by starting with $f_{1}(0)=\cdots=f_{1}(n)=1$
and computing $f_{i}(0),\ldots,f_{i}(n)$ in turn for $i=2,\ldots,d$
(Fig.~\ref{fig: dependency graphs} left). The cost of computing
the values of $f_{i}$ given the previously computed values of $f_{i-1}$
is $O(1+2+\cdots+n)=O(n^{2})$. Thus the complexity of evaluating
$f_{d}(n)$ using (\ref{eq: Richmond}) is $O(n^{2}d)$, a huge improvement
over (\ref{eq: basic formula for AS count}) when $n$ and $d$ are
both small. In contexts where $d$ is very large, however, (\ref{eq: basic formula for AS count})
is impractical. 

\section{An Alternative Recursive Formula}

\begin{figure}
\begin{centering}
\includegraphics{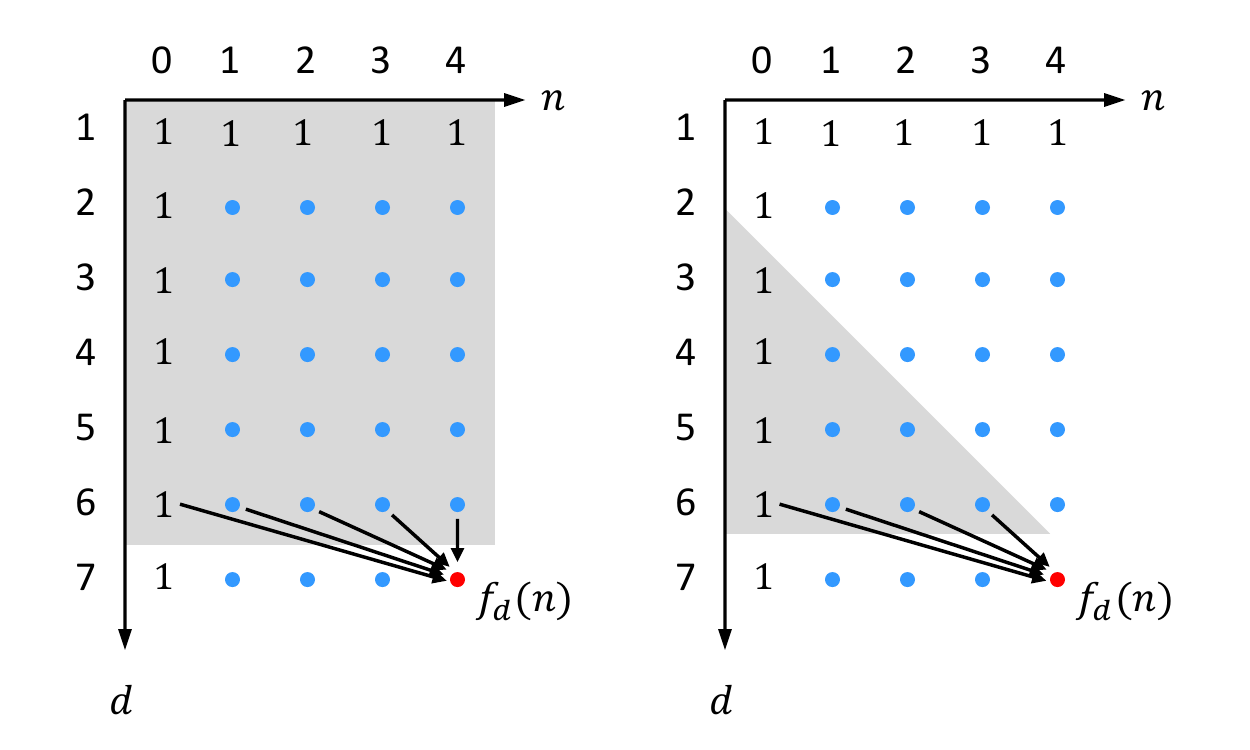}
\par\end{centering}
\caption{\label{fig: dependency graphs}Computational dependencies for two
different recursive formulas for $f_{d}(n)$, the number of abelian
squares. (left) Dependency graph for eq.~(\ref{eq: Richmond}), obtained
from \cite{Richmond2009}. (right) Dependency graph for eq.\,(\ref{eq: new result}).
In each case, the desired quantity $f_{d}(n)$ is shown as a red dot,
arrows show the direct dependencies, and the gray shaded region covers
all the quantities that must be calculated to determine $f_{d}(n)$.
The pattern on the left leads to a cost of $O(n^{2}d)$, while the
pattern on the right leads to a cost of $O(n^{2}\min(n,d))$.}
\end{figure}

In this section I derive an alternative to (\ref{eq: Richmond}) whose
cost of evaluation is only $O(n^{2}\min(n,d))$. Let $A_{d}$ denote
an alphabet of $d$ symbols. The number of abelian squares $(x,y)\in A_{d}^{n}\times A_{d}^{n}$
can be expressed as the sum of the number of anagrams of each word
$x$:
\begin{align}
f_{d}(n) & =\sum_{x\in A_{d}^{n}}\binom{n}{m_{1},\ldots,m_{d}}.
\end{align}
Here $m$ implicitly denotes the signature of $x=(x_{1},\ldots,x_{n})$.
We split off the sum over $x_{n}$:
\begin{align}
f_{d}(n) & =\sum_{x^{\prime}\in A_{d}^{n-1}}\sum_{x_{n}\in A_{d}}\binom{n}{m_{1}^{\prime},\ldots,m_{x_{n}}^{\prime}+1,\ldots m_{d}^{\prime}}
\end{align}
where $m^{\prime}$ is the signature of $x^{\prime}\equiv(x_{1},\ldots,x_{n-1})$.
We have
\begin{align}
\binom{n}{m_{1}^{\prime},\ldots,m_{x_{n}}^{\prime}+1,\ldots m_{d}^{\prime}} & =\frac{n}{m_{x_{n}}^{\prime}+1}\binom{n-1}{m_{1}^{\prime},\ldots,m_{d}^{\prime}}.
\end{align}
Then
\begin{align}
f_{d}(n) & =\sum_{x^{\prime}\in A_{d}^{n-1}}\sum_{x_{n}\in A_{d}}\frac{n}{m_{x_{n}}^{\prime}+1}\binom{n-1}{m_{1}^{\prime},\ldots,m_{d}^{\prime}}.
\end{align}
By symmetry $x_{n}$ can be replaced by any value; choosing $d$ yields
\begin{align}
f_{d}(n) & =d\sum_{x^{\prime}\in A_{d}^{n-1}}\frac{n}{m_{d}^{\prime}+1}\binom{n-1}{m_{1}^{\prime},\ldots,m_{d}^{\prime}}.
\end{align}
Now, each $x^{\prime}$ with a given signature contributes the same
value to the sum. We may thus replace the sum over $x^{\prime}$ by
a sum over the signatures of $x^{\prime}$, weighted by the number
of occurrences of each signature:
\begin{align}
f_{d}(n) & =d\sum_{m_{1}^{\prime}+\cdots+m_{d}^{\prime}=n-1}\frac{n}{m_{d}^{\prime}+1}\binom{n-1}{m_{1}^{\prime},\ldots,m_{d}^{\prime}}^{2}.
\end{align}
We henceforth suppress the primes on $m$. The goal now is to move
the dependence on $m_{d}$ out of the sum, leaving something which
has the form of (\ref{eq: basic formula for AS count}). We have
\begin{align}
\binom{n-1}{m_{1},\ldots,m_{d}} & =\binom{n-1}{m_{d}}\binom{n-1-m_{d}}{m_{1},\ldots,m_{d-1}}.
\end{align}
This yields
\begin{align}
f_{d}(n) & =d\sum_{m_{1}+\cdots+m_{d}=n-1}\frac{n}{m_{d}+1}\binom{n-1}{m_{d}}^{2}\binom{n-1-m_{d}}{m_{1},\ldots,m_{d-1}}^{2}.\\
 & =d\sum_{m_{d}=0}^{n-1}\frac{n}{m_{d}+1}\binom{n-1}{m_{d}}^{2}\sum_{m_{1}+\cdots+m_{d-1}=n-1-m_{d}}\binom{n-1-m_{d}}{m_{1},\ldots,m_{d}}^{2}.
\end{align}
In terms of $k\equiv n-1-m_{d}$,
\begin{align}
f_{d}(n) & =d\sum_{k=0}^{n-1}\frac{n}{n-k}\binom{n-1}{n-1-k}^{2}\sum_{m_{1}+\cdots+m_{d-1}=k}\binom{k}{m_{1},\ldots,m_{d}}^{2}.
\end{align}
Comparison of the latter sum to (\ref{eq: basic formula for AS count})
reveals that it is none other than $f_{d-1}(k)$. The remaining quantities
can be simplified as follows:
\begin{align}
\binom{n-1}{n-1-k} & =\binom{n-1}{k},
\end{align}
\begin{align}
\frac{n}{n-k}\binom{n-1}{n-1-k} & =\binom{n}{k}.
\end{align}
Making these substitutions yields the main result:
\begin{align}
f_{d}(n) & =d\sum_{k=0}^{n-1}\binom{n}{k}\binom{n-1}{k}f_{d-1}(k).\label{eq: new result}
\end{align}
Note the close similarity between (\ref{eq: new result}) and (\ref{eq: Richmond}).
The crucial difference is that in (\ref{eq: new result}) the sum
goes up to only $n-1$; that is, each level of recursion decreases
\emph{both} $n$ and $d$ (Fig.~\ref{fig: dependency graphs} right).
Thus only $\min(n,d)$ levels of recursion are needed. The cost of
this algorithm is $O(n^{2}\min(n,d))$.

Eq. (\ref{eq: new result}) can be interpreted in terms of the following
approach approach to constructing an abelian pair: There are $d$
choices for the first symbol $a$ of $x$. Let $k\in\{0,\ldots,n-1$\}
be the number of occurrences in each word of symbols from $A_{d}/a$.
There are $\tbinom{n-1}{k}$ choices to place those other symbols
in $x$ and $\tbinom{n}{k}$ places to place those other symbols in
$y$. Then, one creates an abelian pair of size $(k,k)$ over $A_{d}/a$,
which is an alphabet of size $d-1$.

\begin{figure}
\begin{centering}
\includegraphics{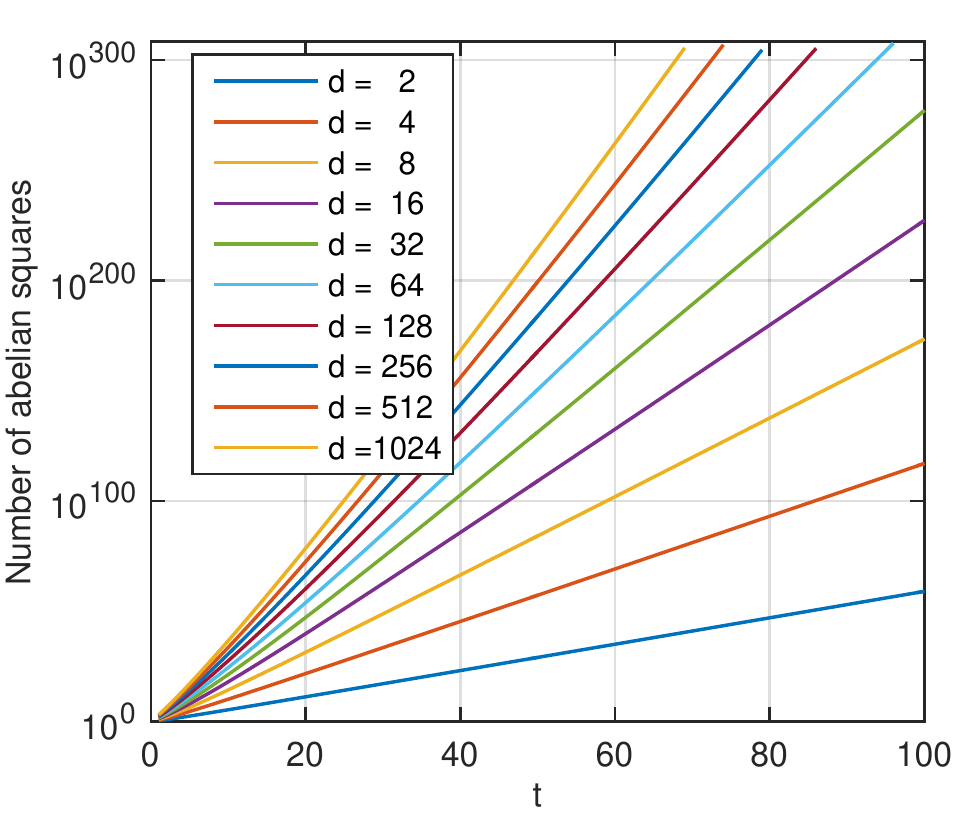}
\end{centering}
\caption{\label{fig: expressivness plots}Number $f_{d}(n)$ of abelian squares
of length $n+n$ over an alphabet of size $d$.}

\end{figure}

Fig. \ref{fig: expressivness plots} shows$f_{d}(n)$ as a function
of $n$ for exponentially increasing values of $d$. (The lines for
$d\ge64$ are truncated due to some of the results being outside the
range of double-precision arithmetic.) The entire plot, comprising
1000 computed points, took less than two seconds to compute in MATLAB
on a standard computer.


\section{Acknowledgements}

This work was performed at Oak Ridge National Laboratory, operated
by UT-Battelle, LLC under contract DE-AC05-00OR22725 for the US Department
of Energy (DOE). Support for the work came from the DOE Advanced Scientific
Computing Research (ASCR) Accelerated Research in Quantum Computing
(ARQC) Program under field work proposal ERKJ354.

\bibliographystyle{unsrtnat}

\input{counting_abelian_squares_biblio.bbl}
\end{document}

%% file: macros.tex
\global\long\def\sinc{\operatorname{sinc}}%

\global\long\def\sech{\operatorname{sech}}%

\global\long\def\rank{\operatorname{rank}}%

\global\long\def\nullsp{\operatorname{null}}%

\global\long\def\vspan{\operatorname{span}}%

\global\long\def\Tr{\operatorname{Tr}}%

\global\long\def\Re{\operatorname{Re}}%

\global\long\def\Im{\operatorname{Im}}%

\global\long\def\ket#1{|#1\rangle}%

\global\long\def\bra#1{\langle#1|}%

\global\long\def\braket#1#2{\langle#1|#2\rangle}%

\global\long\def\Ket#1{\left|#1\right\rangle }%

\global\long\def\Bra#1{\left\langle #1\right|}%

\global\long\def\cket#1{|#1]}%

\global\long\def\cbra#1{[#1|}%

\global\long\def\imag{\mathrm{i}}%

\global\long\def\deriv{\mathrm{d}}%

{\LARGE{}}
\global\long\def\op#1{\operatorname\{#1\}}%
{\LARGE\par}

{\LARGE{}}
\global\long\def\Var{\operatorname{Var}}%
{\LARGE\par}

{\LARGE{}}
\global\long\def\Vol{\operatorname{Vol}}%
{\LARGE\par}

{\LARGE{}}
\global\long\def\Cov{\operatorname{Cov}}%
{\LARGE\par}

{\LARGE{}}
\global\long\def\sign{\operatorname{sign}}%
{\LARGE\par}

{\LARGE{}}
\global\long\def\eig{\operatorname{eig}}%
{\LARGE\par}

{\LARGE{}}
\global\long\def\diag{\operatorname{diag}}%
{\LARGE\par}

\global\long\def\wt{\operatorname{wt}}%

{\LARGE{}}
\global\long\def\Nb{\operatorname{Nb}}%
{\LARGE\par}

\global\long\def\ch{\operatorname{child}}%

\global\long\def\parent{\operatorname{par}}%

{\LARGE{}}
\global\long\def\Perm{\operatorname{Perm}}%
{\LARGE\par}

{\LARGE{}}
\global\long\def\perm{\operatorname{perm}}%
{\LARGE\par}

\global\long\def\cfg{\operatorname{cfg}}%

{\LARGE{}}
\global\long\def\urtri{\urtriangle}%
{\LARGE\par}

%% file: counting_abelian_squares.bbl
\begin{thebibliography}{9}
\providecommand{\natexlab}[1]{#1}
\providecommand{\url}[1]{\texttt{#1}}
\expandafter\ifx\csname urlstyle\endcsname\relax
  \providecommand{\doi}[1]{doi: #1}\else
  \providecommand{\doi}{doi: \begingroup \urlstyle{rm}\Url}\fi

\bibitem[Erd{\H o}s(1957)]{Erdos1957}
Paul Erd{\H o}s.
\newblock Some unsolved problems.
\newblock \emph{Michigan Math. J.}, 4\penalty0 (3):\penalty0 291--300, 1957.
\newblock ISSN 0026-2285, 1945-2365.
\newblock \doi{10.1307/mmj/1028997963}.

\bibitem[Ker{\"a}nen(1992)]{Keranen1992}
Veikko Ker{\"a}nen.
\newblock Abelian squares are avoidable on 4 letters.
\newblock In W.~Kuich, editor, \emph{Automata, {{Languages}} and
  {{Programming}}}, Lecture {{Notes}} in {{Computer Science}}, pages 41--52,
  {Berlin, Heidelberg}, 1992. {Springer}.
\newblock ISBN 978-3-540-47278-0.
\newblock \doi{10.1007/3-540-55719-9-62}.

\bibitem[Iliopoulos et~al.(1997)Iliopoulos, Moore, and Smyth]{Iliopoulos1997}
Costas~S. Iliopoulos, Dennis Moore, and W.~F. Smyth.
\newblock A characterization of the squares in a {{Fibonacci}} string.
\newblock \emph{Theoretical Computer Science}, 172\penalty0 (1):\penalty0
  281--291, February 1997.
\newblock ISSN 0304-3975.
\newblock \doi{10.1016/S0304-3975(96)00141-7}.

\bibitem[Carpi(1998)]{Carpi1998}
Arturo Carpi.
\newblock On the number of {{Abelian}} square-free words on four letters.
\newblock \emph{Discrete Applied Mathematics}, 81\penalty0 (1):\penalty0
  155--167, January 1998.
\newblock ISSN 0166-218X.
\newblock \doi{10.1016/S0166-218X(97)88002-X}.

\bibitem[Cassaigne et~al.(2011)Cassaigne, Richomme, Saari, and
  Zamboni]{Cassaigne2011}
Julien Cassaigne, Gw{\'e}na{\"e}l Richomme, Kalle Saari, and Luca~Q. Zamboni.
\newblock Avoiding abelian powers in binary words with bounded abelian
  complexity.
\newblock \emph{Int. J. Found. Comput. Sci.}, 22\penalty0 (04):\penalty0
  905--920, June 2011.
\newblock ISSN 0129-0541.
\newblock \doi{10.1142/S0129054111008489}.

\bibitem[Huova et~al.(2012)Huova, Karhum{\"a}ki, and Saarela]{Huova2012}
Mari Huova, Juhani Karhum{\"a}ki, and Aleksi Saarela.
\newblock Problems in between words and abelian words: K-abelian avoidability.
\newblock \emph{Theoretical Computer Science}, 454:\penalty0 172--177, October
  2012.
\newblock ISSN 0304-3975.
\newblock \doi{10.1016/j.tcs.2012.03.010}.

\bibitem[Crochemore et~al.(2013)Crochemore, Iliopoulos, Kociumaka, Kubica,
  Pachocki, Radoszewski, Rytter, Tyczy{\'n}ski, and Wale{\'n}]{Crochemore2013}
M.~Crochemore, C.~S. Iliopoulos, T.~Kociumaka, M.~Kubica, J.~Pachocki,
  J.~Radoszewski, W.~Rytter, W.~Tyczy{\'n}ski, and T.~Wale{\'n}.
\newblock A note on efficient computation of all {{Abelian}} periods in a
  string.
\newblock \emph{Information Processing Letters}, 113\penalty0 (3):\penalty0
  74--77, February 2013.
\newblock ISSN 0020-0190.
\newblock \doi{10.1016/j.ipl.2012.11.001}.

\bibitem[Bennink((submitted))]{Bennink2022b}
Ryan Bennink.
\newblock Counting abelian squares for a problem in quantum computing.
\newblock \emph{Journal of Combinatorics}, (submitted).

\bibitem[Richmond and Shallit(2009)]{Richmond2009}
L.~B. Richmond and Jeffrey Shallit.
\newblock Counting {{Abelian Squares}}.
\newblock \emph{Electron. J. Combin.}, 16\penalty0 (1):\penalty0 R72, June
  2009.
\newblock ISSN 1077-8926.
\newblock \doi{10.37236/161}.

\end{thebibliography}
